\documentclass[a4paper, 11pt]{article}
\usepackage{graphicx}

\usepackage{wrapfig}
\usepackage{floatflt}
\usepackage[english]{babel}
\usepackage[usenames]{color}

\usepackage{amsfonts}
\usepackage{amsmath}
\usepackage{amssymb}
\usepackage{mathtext}

\newcounter{theorem}
\newenvironment{theorem}
{\mbox{} \newline   \refstepcounter{theorem} {\bf Theorem \thetheorem.} \em   }
{}
\renewcommand{\thetheorem}{\arabic{theorem}}

\newcounter{definition}
\newenvironment{definition}
{\mbox{} \newline   \refstepcounter{definition} {\bf Definition \thedefinition.} \em   }
{}
\renewcommand{\thedefinition}{\arabic{definition}}

\newcounter{statement}

\renewcommand{\thestatement}{\arabic{statement}}

\newcounter{lemma}

\renewcommand{\thelemma}{\arabic{lemma}}

\newcounter{supposition}

\renewcommand{\thesupposition}{\arabic{supposition}}

\newcounter{consequence}

\renewcommand{\theconsequence}{\arabic{consequence}}

\newcounter{remark}

\newcommand{\eps}{\varepsilon}

\newcommand{\tr}{\mathop{\rm tr}}

\newcommand{\area}{\mathop{\rm area}}

\renewcommand{\geq}{\geqslant}
\renewcommand{\leq}{\leqslant}

\newcommand{\mystar}{\mathop{\rm star}}

\newcommand{\tg}{\mathop{\rm tg}}

\newcommand{\R}{\mathbb{R}}

\def\p{\partial}

\newenvironment{myfigure}[3]
{ \vbox{ \refstepcounter{figure} \label{#1} \begin{center} {#2} \end{center}
\begin{center}  {\footnotesize Fig.~\ref{#1}. {#3} }\end{center} \vspace{-1ex} }
{       }                                                     }

\title{\bf Discrete extrinsic curvatures based on polar polyhedra concept}
\author{V.\,V.\,Garanzha\footnotemark[1] \\[2mm]
{Computing Center RAS, garan@ccas.ru}}

\date{}

\begin{document}

\maketitle

\footnotetext[1]{Research supported by grant OMN-03 of Department of mathematical sciences, Russian academy of sciences and by program "Leading Scientific
Schools" (project no. NSh-5073.2008.1)}

\begin{abstract}
Duality principle for approximation of geometrical objects (also known as Eudoxus exhaustion
method) was extended and perfected by Archimedes in his famous tractate ``Measurement of circle''.
The main idea of the approximation method by Archimedes is to construct a sequence of pairs of
inscribed and circumscribed polygons (polyhedra) which approximate curvilinear convex body. This
sequence allows to approximate length of curve, as well as area and volume of the bodies and to
obtain error estimates for approximation. In this work it is shown that a sequence of pairs of
locally polar polyhedra allows to construct piecewise-affine approximation to scherical Gauss map,
to construct convergent pointwise approximations to mean and Gauss curvature, as well as to obtain
natural discretizations of bending energies.
\end{abstract}


{\bf Discrete curvature functionals and surfaces of bounded curvature. }  One of the hard problems
of modern geometry is approximation of nonregular surfaces by polyhedra. In the sense of intrinsic
metric (based on distance along surface) this problem was solved in the works of A.D. Alexandrov
and his scientific school \cite{Alexandrov-1962}. A.D. Alexandrov developed theory of ``good''
approximation of manifolds of bounded curvature by polyhedral manifolds. However these results are
not sufficient to establish ``good'' convergence in the sense of extrinsic metric. The class of
surfaces being manifolds of bounded curvatures in the intrinsic sense is well defined: they are
called surfaces of bounded curvature\cite{Alexandrov-1950-3c},  \cite{Bakelman-1956-6a},
\cite{Pogorelov-1956i}, \cite{Burago-1968}.

Extrinsic curvatures for polyhedra are introduced using integral relations. Gauss-Bonnet  theorem
allows to assign to the vertex of polyhedron  curvature which is equal to angular excess of its
conical neighborhood \cite{Alexandrov-1950}, \cite{Alexandrov-1962}. Balance equations for  vector
mean curvature can be used to derive discrete mean curvatures for polyhedra and to construct
discrete approximation to Laplace-Beltrami operator.  To this end one can also use variation of
surface area and its relation with the sweep volume  \cite{Sullivan-2002}. In
\cite{Cohen-Steiner-2003}  with each region on the surface it is associated a tensor which in the
smooth case is the average of curvature tensor over this region. For polyhedral domain the same
value provides weakly convergent estimator of the curvature tensor.  In \cite{Hildebrandt-2006} it
is shown that if a sequence of polyhedral surfaces converges to a regular surface in Hausdorff
distance, then the following conditions are equivalent: a) convergence of normal fields, b)
convergence of metric tensors, c) convergence of area, d) convergence of Laplace-Beltrami
operators.

In order to define ``good'' approximation by polyhedral surface, one have to define such curvature
measures for non-regular surface $M$, which can be introduced via sequence of polyhedral surfaces
$P_k$ which the following properties: a) $P_k$ converge to $M$ pointwise when $k \to + \infty$; b)
$P_k$ converge to $M$ uniformly in intrinsic metric; b) positive and negative parts of curvature
of $P_k$ converge to to positive and negative part of curvature of $M$ in a weak sense; c)
spherical Gauss map of $P_k$ converge to spherical map of $M$ in a weak sense. Approximation
problem in such a setting is still not solved \cite{discrete-geometry-2008}.

Consider curvature functionals which can be used for investigation of non-regular surfaces. Let
$M$ denote regular 2D surface in $\R^3$ (regular in a sense that it admits thrice continuously
differentiable local parameterization). Consider functional
\begin{equation} \label{eq.curv-generic-curvature-measure}
E_g(M)  = \frac12 \int\limits_M  g(A) \, d \sigma, \ \ g(A) \geq |\det A|,
\end{equation}
where  $d \sigma$ is the surface area element, $A \in \R^{2 \times 2}$ is the matrix of the shape operator or curvature tensor defined by equality
\[
A =  G^{-1}T,
\]
Here $G$ and $T$ are matrices of the first and second fundamental forms of the surface, respectively, and $g(A)$ is certain curvature density measure. Functional $E_g(M)$ is called bending energy of the surface. If energy $E_g(M)$ is bounded, then absolute Gauss curvature of the surface is bounded as well.

Matrix $A$ is nothing else but the jacobian matrix of the spherical map $\mu$. Let us remind that spherical map identifies with each point $p$ of regular surface $M$ a point $b = \nu(p)$ on a unit sphere $\mathbb S^2$, where $\nu(p)$ is the unit normal to surface. One can compute intersection point $q$ between plane, passing through $b$ and orthogonal to $\nu(p)$ and a ray going through $\nu(p^\prime)$, where point $p^\prime$ belongs to some neighborhood of $p$,  which is shown in fig.~\ref{fig.curv-normal-graph} a). Mapping $p^\prime \to q$ is called normal map and defines normal image of neighborhood  of $p$.


\begin{myfigure}{fig.curv-normal-graph}
{
 a) \hspace{4cm} b) \\
 \mbox{\includegraphics[scale=0.5, angle=0]{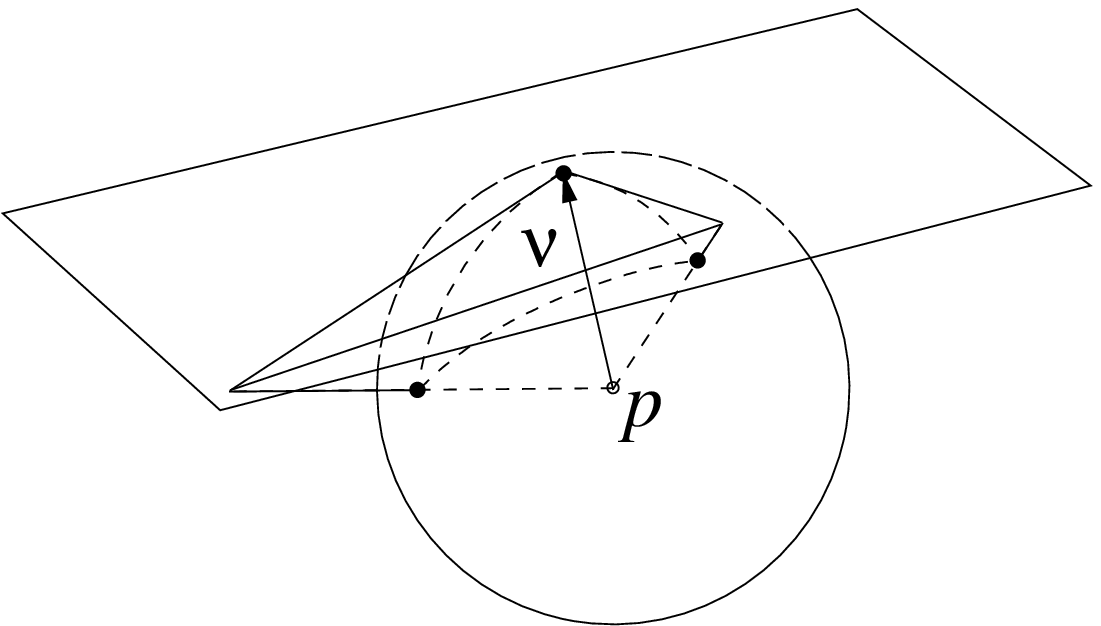} \includegraphics[scale=0.4, angle=0]{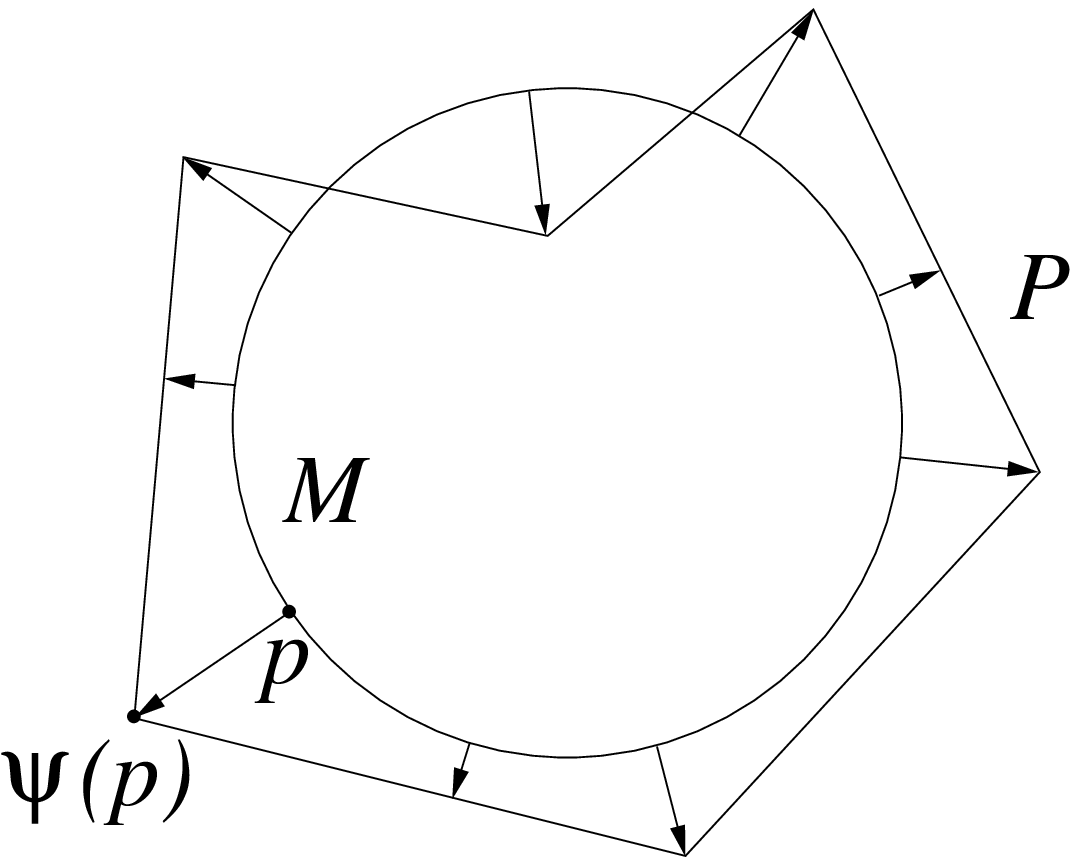}}
}{a) Spherical and normal images, b) normal graph over the surface. }
\end{myfigure}

It is said that polyhedral surface $P_h$ is normal graph over $M$ if projection $\psi$ of surface $M$ onto $P_h$ along normals to $M$ is homemorphism, which is shown on fig.~\ref{fig.curv-normal-graph} b).

Well known example of bending energy is given by mean total quadratic curvature measure
\[
 E_2(M) =  \frac12 \int\limits_M  \tr (A^T A) \, d \sigma
\]
Absolute minimum of this functional is attained when a surface homeomorphic to sphere is precisely the sphere. Mean quadratic curvature measure is not suitable for description of nonregular surfaces since it is not defined for polyhedra. In other words it takes infinite value for polyhedral surface.

If bending energy  majorates absolute curvature and remains bounded for refined sequence of polyhedra then one can expect that the limiting surface for this sequence will be surface of bounded curvature.

One can consider the following curvature measures which make sense for polyhedra:
\begin{equation} \label{eq.curv-mean+absolute}
E_1(M) =  \int\limits_M (\left( \tr(A^T A) \right)^\frac12 + |\det A|)  \, d \sigma
\end{equation}
and
\begin{equation} \label{eq.curv-combined-energy}
E_\eps(M) =  \frac1{\sqrt{2}} \int\limits_M  \tr(A^T A) \frac{(\eps + |\det A|)^\frac12}{(\eps + \tr(A^T A))^\frac12 }  \, d \sigma,
\end{equation}
where $\eps > 0$ is a constant.

\noindent
{\bf Duality principle and approximation of surfaces by polyhedra.} Let us consider the method for construction of discrete bending energies suitable for approximation of nonregular surfaces.

Consider 2D paraboloid
\[
P = \{x: \  x_3 = u(x_1, x_2), u(x_1, x_2) =
\]
\[
= \frac12 (h_{11} x_1^2 + 2 h_{12} x_1 x_2 + h_{22} x_2^2) \}
\]
We shall use upper index $l$ to denote vectors from $\R^3$, while values without $l$ will denote their orthogonal projections onto plane $x_3=0$. It is convenient to write function $u$ as  $u(p) = \frac12 p^T H p$, where $H$ is the shape operator matrix of paraboloid $P$ at the origin.


\begin{myfigure}{fig.curv-polar-domain}
{
a) \hspace{4cm} b) \\
\includegraphics[scale=0.4, angle=0]{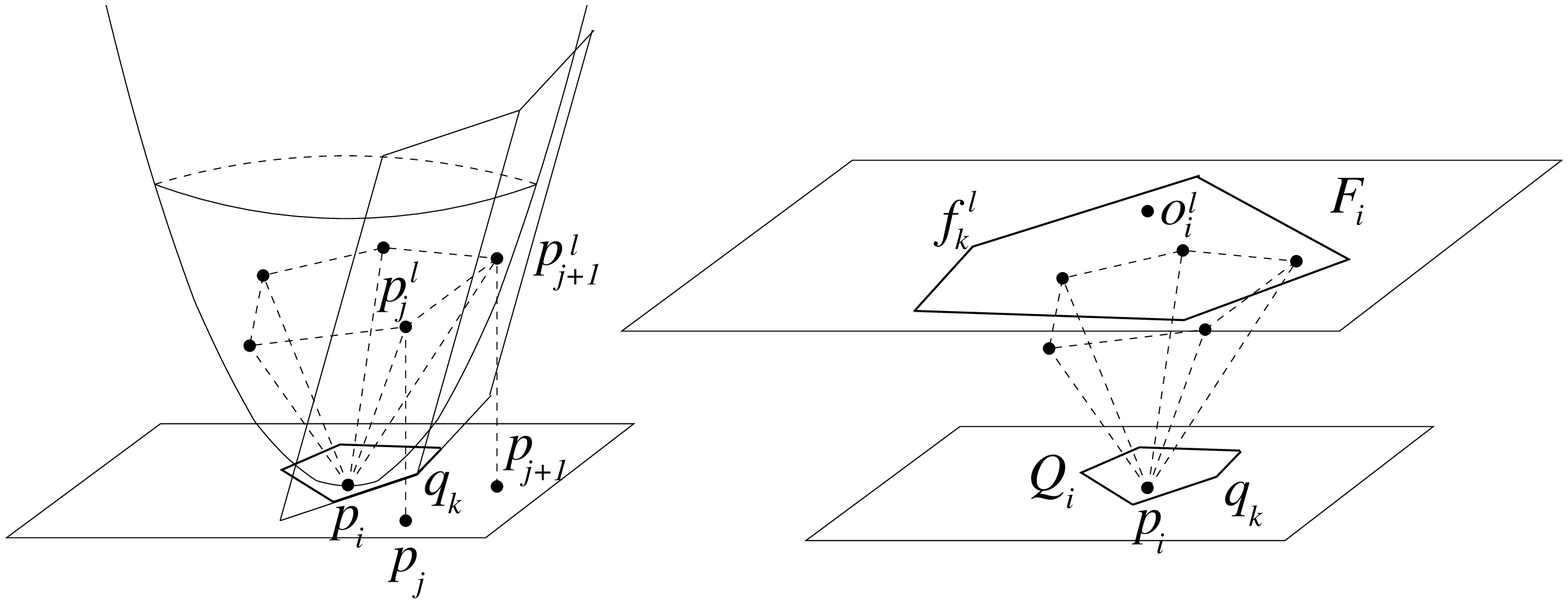}
}{ a) Polyhedral surface inscribed into elliptic paraboloid; b) dual face and normal image of the vertex.}
\end{myfigure}

Consider fragment of convex polyhedral surface $P_h$ inscribed into elliptic paraboloid. The faces of this polyhedron incident to vertex $p^l_i=0$ are shown in fig.~\ref{fig.curv-polar-domain} a). The plane $x_3=0$ is tangent plane at the origin. Tangent planes passing through the vertices $p_j^l$ lying at the edges, incident to $p^l_i$, cut on the plane $x_3=0$ polygon $Q_i$. In fact $Q_i$ is the face of dual (polar) polyhedra $P_h^\star$ circumscribed around the same paraboloid. Let us denote by $\mathcal V(p^l_i)$ the set of vertices of $P_h$ belonging to edges, incident to  $p^l_i$, while notation $\mathcal V(G)$ is used for the set of vertices of the face $G$.

One can construct normal image of the vertex $p^l_i$, namely the convex polygon $F_i$ on the plane $x_3=1$. The vertices of this polygon are intersections of rays passing through $p^l_i$ and orthogonal to faces, incident to $p^l_i$, with the plane $x_3=1$. Polygons $Q_i$ and $F_i$ are shown in fig.~\ref{fig.curv-polar-domain} b).

Consider vertex $q_k^l$ of the polygon $Q_i$. Vector $q^l_k$  can be found as the solution of the linear system
\begin{equation} \label{eq.curv-vertex-dual-face}
\begin{array}{l}
{n^l(p^l_i)}^T (q^l_k - p^l_i) = 0 \\
{n^l(p^l_j)}^T (q^l_k - p^l_j) = 0, \ j \in \mathcal V(G_k), \ i \neq j
\end{array}
\end{equation}
where $n^l$ are unscaled normals to paraboloid and $p^l_j$ are vertices of a face $G_k \in \mystar(p^l_i)$, i.e.
incident to $p^l_i$.

Linear system (\ref{eq.curv-vertex-dual-face}) is nonsingular. In order to obtain normal $n^l$,  one can compute gradient of the function $x_3 - u(x_1, x_2)=0$, i.e.,
\[
n^l(p_j^l) =
\left(
\begin{array}{l}
-H p_j \\ 1
\end{array}
\right),  \ \
n^l(p^l_i) =
\left(
\begin{array}{l}
0 \\ 1
\end{array}
 \right)
\]
Thus linear system (\ref{eq.curv-vertex-dual-face}) can be reduced to
\[
( -Hp_j)^T (q_k - p_j) - u(p_j) = 0,
\]
or
\begin{equation} \label{eq.curv-q-equation}
(Hq_k)^T p_j = u(p_j), \  j \in \mathcal V(G_k), \ i \neq j
\end{equation}
It can be easily verified that if the number of equations in system (\ref{eq.curv-q-equation}) exceeds the number of unknowns it remains consistent.

Now let us find vertex $f^l_k$ of the polygon $F_i$. Its coordinates are the solution of linear system
\begin{equation} \label{eq.curv-normale-image-vertex}
\begin{array}{c}
n^l(p^l_i)^T (f^l_k - p^l_i) = |n^l(p^l_i)|, \\
(p^l_j - p^l_i)^T (f^l_k - p^l_i) = 0, \ j \in \mathcal V(G_k), \ i \neq j
\end{array}
\end{equation}
This system simplifies to
\begin{equation} \label{eq.curv-f-equation}
p_j^T f_k  + u(p_j) = 0,  \  j \in \mathcal V(G_k), \ i \neq j
\end{equation}
From equations (\ref{eq.curv-q-equation}),  (\ref{eq.curv-f-equation}) follows equality
\begin{equation} \label{eq.curv-f=Hq}
f_k = -H q_k,
\end{equation}
which hold for all vertices of polygon $F_i$. Thus the following theorem is proven:
\begin{theorem} (\cite{Garanzha-2008})  \label{th.curv-affine-map-curvature}
polygons $F_i$ and $Q_i$ are affine equivalent, i.e. $Q_i = \phi_i^\star (F_i)$ and jacobian matrix of affine map $\phi^\star_i$ coincides with $-H$, where $H$ is the matrix of shape operator of paraboloid $P$ at the origin.
\end{theorem}

Let us note that deriving equality $f_k = -Hq_k$ we did not used the fact that matrix $H$ is positive definite. Formally (\ref{eq.curv-f=Hq}) holds for arbitrary matrix $H$. It is just required that matrix with vectors $p_j$ as columns has the full rank. Thus duality principle for computation of curvature tensor can be applied in the case of hyperbolic paraboloid, shown in fig.~\ref{fig.curv-polar-domain-concave}.


\begin{myfigure}{fig.curv-polar-domain-concave}
{
a) \hspace{4cm} b) \\
\includegraphics[scale=0.45, angle=0]{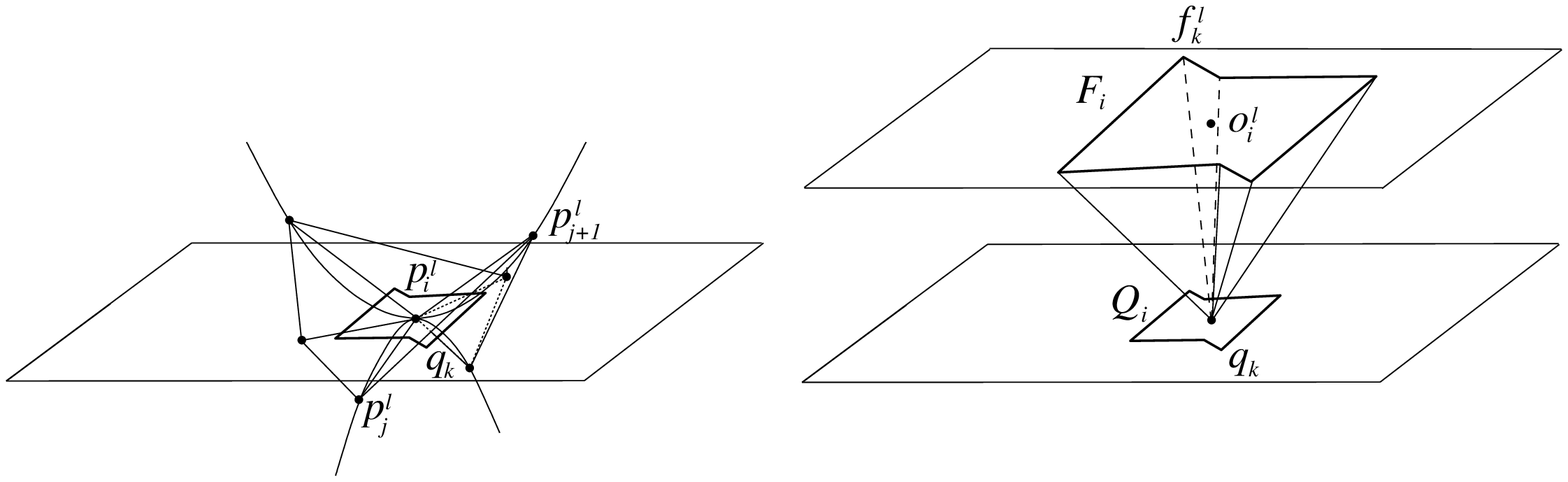}
}{ a) Polyhedral surface inscribed into hyperbolic paraboloid; b) dual face and normal image of the vertex.}
\end{myfigure}

From the duality principle it follows that face $G_k$ of polyhedron $P_h$, incident to $p^l_i$ corresponds to vertices $q_k^l$ and $f_k^l$. If faces $G_m$ and $G_k$ have common edge, then vertices $q_m^l$ and $q_k^l$ should be connected by an edge. The same is true for vertices $f_m^l$ and $f_k^l$. These arguments can be applied in the case when boundary of $F_i$ and $Q_i$ is self-intersecting closed polyline. In this case vertex $p^l_i$ can be called nonregular since polyhedron $P_h$ provides poor local approximation to paraboloid $P$ in the neighborhood of $p^l_i$.

The condition that $P_h$ is inscribed polyhedral surface and $P_h^\star$ is the circumscribed one is the particular case of polarity with respect to paraboloid \cite{Fenchel-1949}. It is well known that relation of polarity, i.e. one-to-one correspondence between point and a plane can be introduced using arbitrary surface $S$  of the second order \cite{Alexandrov-1950}.

Consider rays originating from certain point $p$ and tangent to $S$. All points of contact with S will lie in the same plane which is polar to point $p$. It may happen that point $p$ is situated in such a way that tangent rays cannot originate from it. Then one can draw an arbitrary plane $\Pi$ through $p$. Tangent rays passing through intersection points between $S$ and $\Pi$ consitute a cone. When plane $\Pi$ is varied, the summit of this cone sweep the plane which is precisely the polar plane to point $p$. Let us remark that in convex analysis and optimization only special case of polarity with respect to sphere is considered.

Suppose now that vertices $p^l_j$ of $P_h$ do not lie on the surface of paraboloid $P$, i.e.
\[
(p^l_j)_3 = \frac12 p_j^T H p_j + \delta_j
\]
The plane of the face, polar to $p_j^l$ is defined by equality
\[
(x^l)_3 + (p^l_j)_3 = p_j^T H x
\]
Suppose that $p_i = 0$, then $(p_i^l)_3 = \delta_i$. Point $f_k^l$ is defined by system (\ref{eq.curv-normale-image-vertex}), i.e.
\begin{equation} \label{eq.curv-f-equation-polar}
p_j^T f_k  + \frac12 p_j^T H p_j + \delta_j - \delta_i = 0,  \  j \in \mathcal V(G_k), \ i \neq j
\end{equation}
while $q_k$ is defined as the intersection point of a planes, polar to the vertices of edge $G_k$, thus
\begin{equation} \label{eq.curv-q-equation-polar}
(Hq_k)^T p_j = \frac12 p_j^T H p_j + \delta_j - \delta_i, \  j \in \mathcal V(G_k), \ i \neq j
\end{equation}
As a result one obtains equality $f_k = -H q_k$, which means that theorem~\ref{th.curv-affine-map-curvature} holds for polyhedral surfaces, polar with respect to paraboloid $P$, and not just for inscribed and circumscribed polyhedra.

Now consider the case when there exists face $G_k$ of $P_h$, parallel to the plane $x_3=0$. Vertex $q^l_k$ dual to this face is intersection of planes polar to points $p^l_j$ being vertices of $G_k$, and $q_k = 0$. Consider normal image of the neighborhood of vertex $q^l_k$, i.e. polygon $B_k$ with vertices computed as solutions to linear system
 \begin{equation} \label{eq.curv-normale-image-vertex-dual}
 \begin{array}{c}
n^l(G_k)^T (b^l_j - q^l_k) = |n^l(G_k)|, \\
(q^l_i - q^l_k)^T (b^l_j - q^l_k) = 0, \  i \in \mathcal V(q^l_j), \ i \neq j,
\end{array}
\end{equation}
where $n^l(G_k)$ denotes normal to face $G_k$.
Elementary calculation shows that
\[
b^l_j = - H p^l_j,
\]
i.e. it is proven the following theorem
\begin{theorem} (\cite{Garanzha-2008}) \label{th.curv-affine-map-curvature-dual}
Polygons $B_k$ and $G_k$ are affine equivalent, i.e. $B_k = \phi_k (G_k)$ and jacobian matrix of affine mapping $\phi_k$ coincides with $-H$, where $H$ is the curvature tensor of paraboloid $P$ at the origin.
\end{theorem}

Let us consider arbitrary regular closed 2D surface $M$ and polyhedral surface $P_h$ inscribed into $M$. One can choose cartesian frame $x_i$ such that surface $M$ can be locally written as $x_3 = f(x_1, x_3)$, and
\[
f(x_1, x_2) = u(x_1, x_2) + O(|x|^3), \ \ h_{ij} = \frac{\p^2 f}{\p x_i \p x_j} (0,0),
\]
which means that $P$ is touching paraboloid at the point $p^l_i$ of the surface $M$.
Denote by $F_i$ normal image of neighborhood of the vertex $p^l_i$. The vertices of polygon $F_i$ are computed via solution of linear system (\ref{eq.curv-normale-image-vertex}). Let us construct dual polyhedral surface $P_h^\star$, consisting of the faces tangent to $M$. Tangent face, dual to $p^l_i$ is denoted by $Q_i$. Vertices of polygon $Q_i$ can be found by solving (\ref{eq.curv-vertex-dual-face}), where $n^l$ are normals to the surface $M$. In local coordinate frame they can be written as
\[
n^l(p^l_j) = \left( \begin{array}{c} - \nabla f(p_j)  \\ 1 \end{array}\right)
\]
For a vertex $q^l_k$  of a dual polyhedral surface $P_h^\star$ one can find in turn the normal image, polygon $B_k$. The plane of this polygon is parallel to the face $G_k$, dual to $q^l_k$. Vertices of polygon $B_k$ can be found using equation (\ref{eq.curv-normale-image-vertex-dual}).

\begin{definition} \label{def-regular-1}
The vertex $p^l_i$ of polyhedral surface $P_h$ is called regular if its dual polygon $Q_i$ and orthogonal projection of polygon $F_i$ onto the plane of polygon $Q_i$ are star-shaped domains with respect to the point $p^l_i$.
\end{definition}

\begin{definition} \label{def-regular-2}
The vertex $q^l_k$ of polyhedral surface $P_h^\star$ is called regular if face $G_k$ of $P_h$ and
orthogonal projection of polygon $B_i$ onto the plane of face $G_k$ are star-shaped domains with
respect to the projection of point $q^l_k$.
\end{definition}

It is clear that these two definitions are mutually symmetric.

For non-convex surface the above definitions are too restrictive. One can use definition of weak
regularity:
\begin{definition}
The vertex $p^l_i$ of polyhedral surface $P_h$ is called weakly regular if its dual polygon $Q_i$ and orthogonal projection of polygon $F_i$ onto the plane of polygon $Q_i$ are simple polygons and contain  point $p^l_i$ inside.
\end{definition}
If polygon $B_k$ is simple then one can assume that vertex $q^l_k$  is weakly regular.

\begin{myfigure}{fig.curv-parabola-tri-dual}
{
a) \hspace{3cm} b) \\
\mbox{\includegraphics[scale=0.2, angle=0]{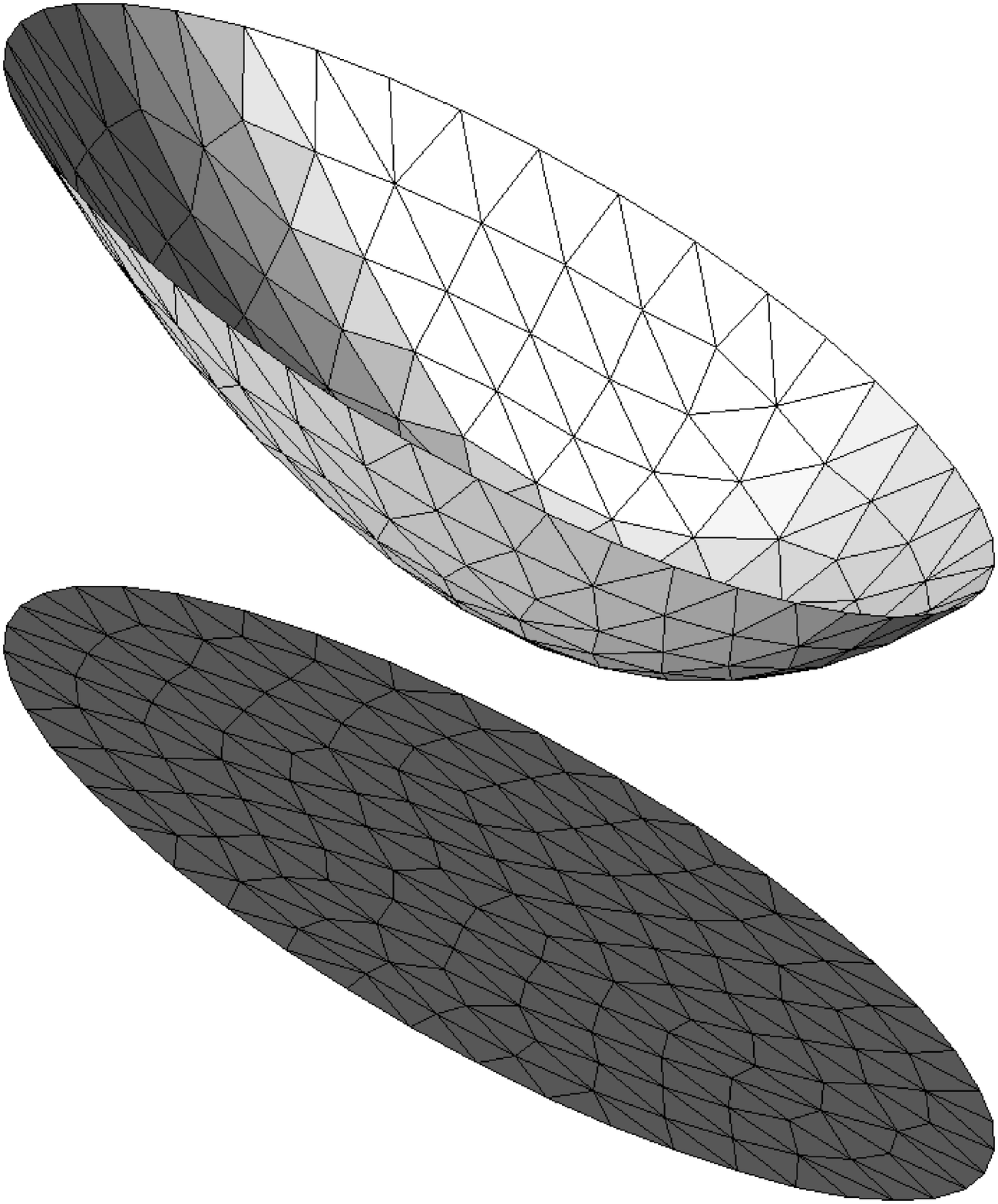} \hspace{-0.25cm} \includegraphics[scale=0.2, angle=0]{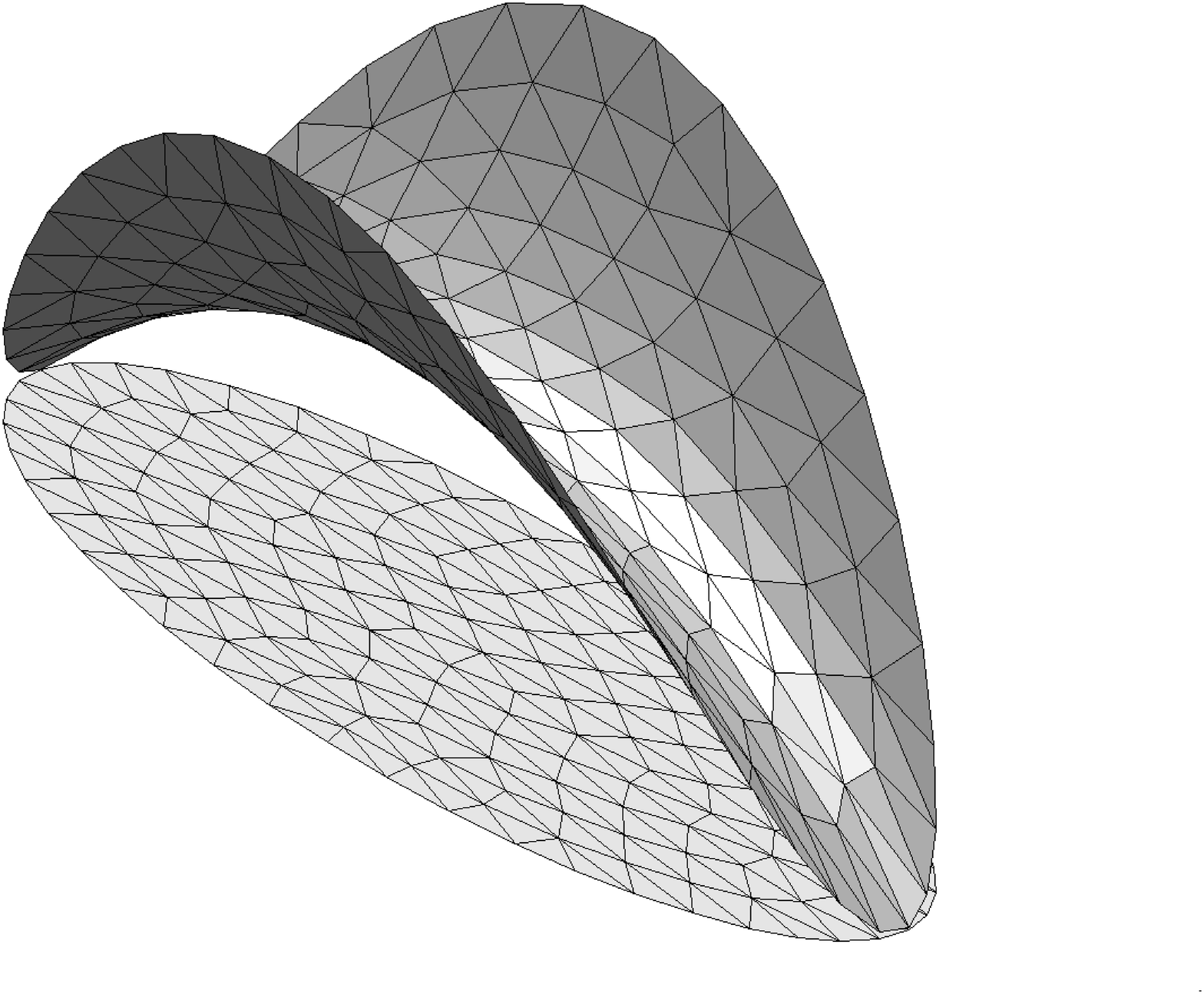}
}
}{Triangulations of elliptic and hyperbolic paraboloids and their projections onto horizontal plane.}
\end{myfigure}

\begin{myfigure}{fig.curv-parabola-tri-dual-1}
{
\mbox{\includegraphics[scale=0.2, angle=0]{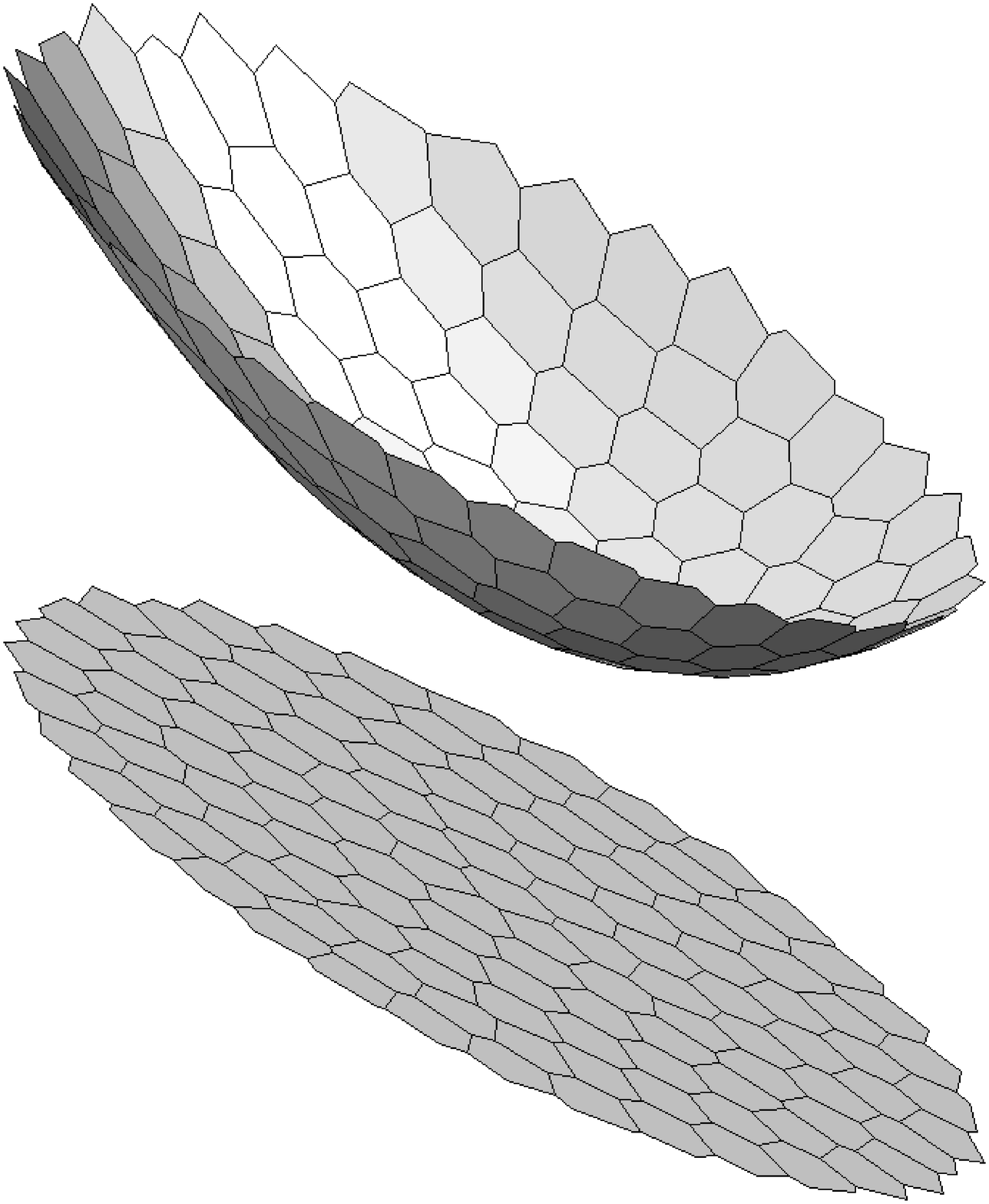}  \hspace{-0.25cm} \includegraphics[scale=0.2, angle=0]{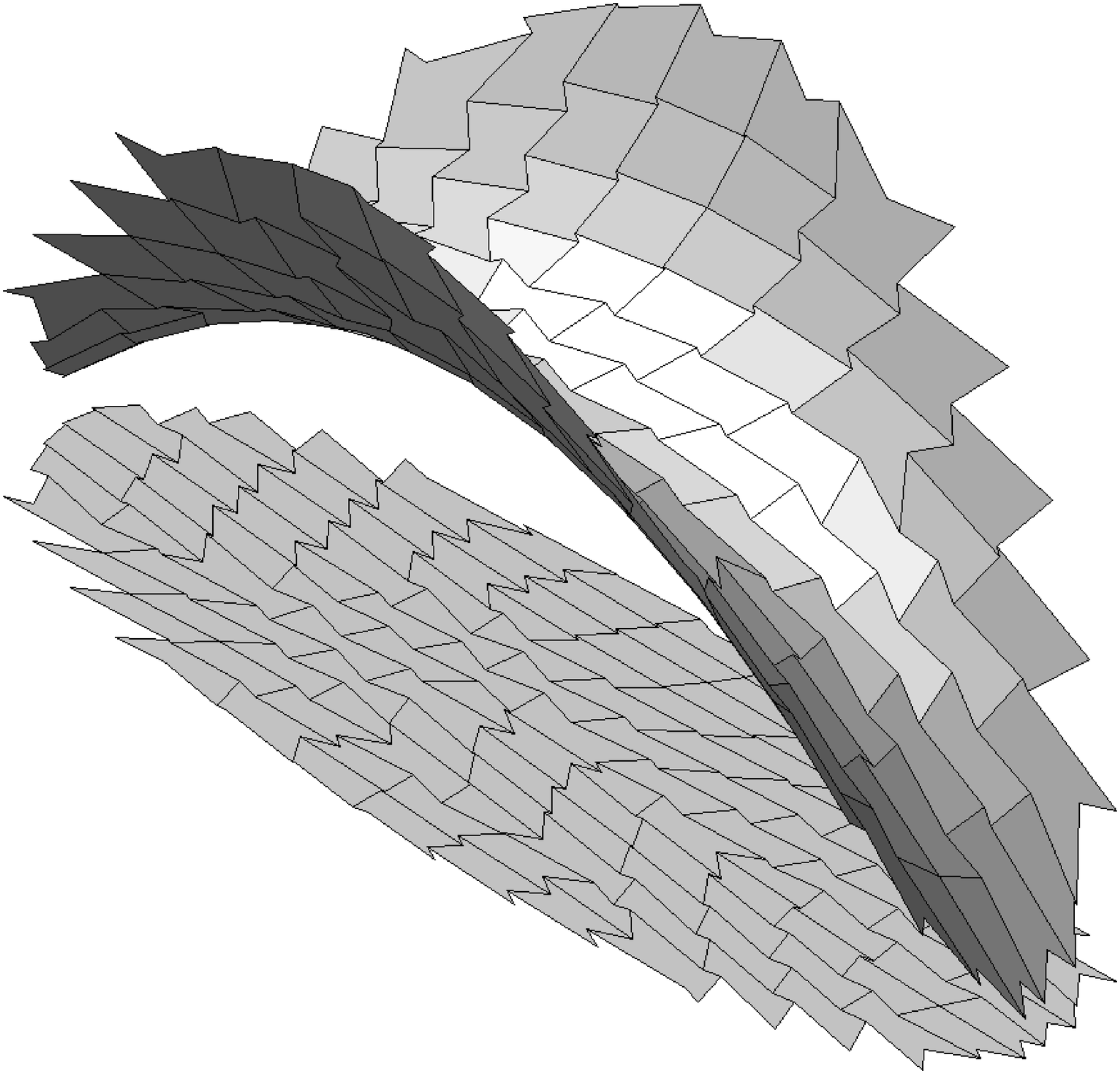}}
}{Dual polyhedral surfaces and their projections.}
\end{myfigure}

On fig.~\ref{fig.curv-parabola-tri-dual}  there are shown fragments of triangulated polyhedral surfaces $P_h$ inscribed into elliptic and hyperbolic paraboloid and their projection on the plane $x_3 = 0$. Dual polyhedral surfaces $P_h^\star$ and their projections are shown in fig.~\ref{fig.curv-parabola-tri-dual-1}. In the case of elliptic paraboloid all dual faces are convex polygons, while in the hyperbolic case dual faces $Q_i$  are quadrilaterals with concave sides(edges), i.e. the turn of edges from the side of $Q_i$ is nonpositive.

It should be noted that dual polyhedral surface $P_h^\star$ in the case of elliptic paraboloid is precisely the Voronoi generatrice and projection of its faces onto horizontal plane make partitioning of a plane into convex polygons being affine image of Voronoi partitioning.

In general case one cannot state that polygons $Q_i$ and $F_i$ are affine equivalent. The polygons $G_k$ and $B_k$ are not affine equivalent as well. Thus one have to construct piecewise affine homemorphisms $\phi_i^\star: Q_i \to F_i$ and $\phi_k: G_k \to B_k$. Without loss of generality one can assume that $\phi_i^\star$ and  $\phi_k$ are 2D mappings. If all vertices of $P_h$ and $P_h^\star$ are regular, then polygons $Q_i$ and $G_k$ can be triangulated simply by connecting their vertices with the point with respect to which they are star-shaped. In the case of weak regularity one have to construct more general triangulations. Let us denote by $\mathcal T^Q_i$ triangulation of $Q_i$, and by $\mathcal T^G_k$ - triangulation of $G_k$. The images of these triangulations under maps $\phi^\star_i$ and $\phi_k$ are triangulations $\mathcal T^F_i$ and  $\mathcal T^B_k$. The number of triangles can be sharply reduced if all faces $G_k$ are triangles. Then mapping $\phi_i: G_k \to B_k$ can be taken as affine one.

Denote by $-A^\star_{im}$ jacobian matrix of affine map of triangle $T^\star_{im} \in \mathcal T^Q_i$ onto $m$-th triangle of $\mathcal T^F_i$, i.e.
\begin{equation} \label{eq.curv-jacobi-matrix-dual}
A^\star_{im} = - \nabla \phi_i^\star |_{T^\star_{im}}
\end{equation}
Denote by $-A_{km}$ jacobian matrix of affine map of triangle $T_{km} \in \mathcal T^G_k$
onto $m$-th triangle of $\mathcal T^B_i$, i.e.
\begin{equation} \label{eq.curv-jacobi-matrix}
A_{km} = - \nabla \phi_k |_{T_{km}}
\end{equation}

It is obvious that if all vertices of $P_h$ and $P_h^\star$ are regular, and diameters of faces of $P_h$ tend to zero, matrices $A^\star_{im}, A_{km}$ converge to exact value of the curvature tensor $A$,  i.e.
\[
A^\star_{im} \to A(p^l_i), \ \   A_{km} \to A({\psi^\star}^{-1}(q^l_k))
\]
where $\psi^\star$ is the homeomorphism which maps $M$ onto $P_h^\star$  along normals to $M$.

Let us remark that matrices $A^\star_{im}$ and $A_{km}$ in general are not symmetric. Hence in order to compute principal curvatures one have to use singular values of these matrices instead of eigenvalues, and use SVD (singular value decomposition) of $A^\star_{im}$ and $A_{km}$ in order to compute principal directions.

In practice exact surface normals are not known. Many methods are available for computation of normals at the vertices of polyhedra. One possible way to find normals approximately is to require that discrete curvature of  $P_h$ is as close as possible to discrete curvature of dual surface $P_h^\star$. In order to make this comparison possible, one has to construct piecewise-affine homeomorphism $\psi_h: P^\star_h \to P_h$.

\begin{myfigure}{fig.curv-integration-domain-convex}
{
a) \hspace{4cm} b)  \\
\mbox{\includegraphics[scale=0.5, angle=0]{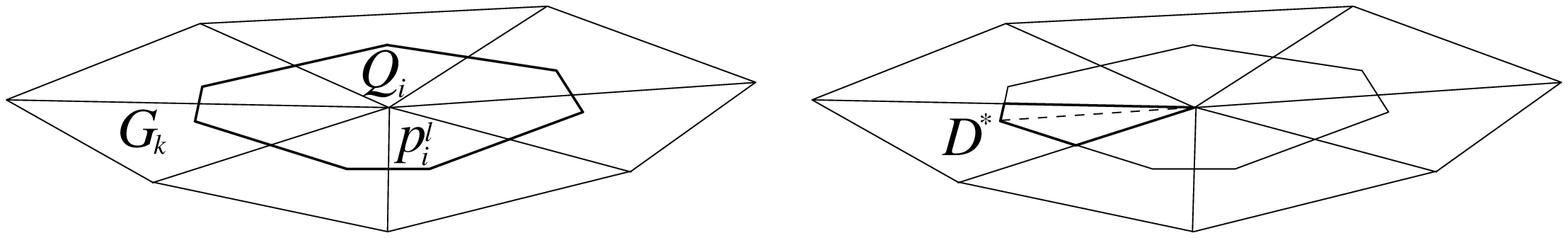}
}
}{Construction of piecewise-affine homeomorphism $\psi_h: P_h^\star \to P_h$.}
\end{myfigure}

On fig.~\ref{fig.curv-integration-domain-convex} a) it is shown fragment of polyhedral surface, with direction of view being orthogonal to the plane of polygon $Q_i$. Projection of the face $G_k$ onto the plane of dual face $Q_i$ has nonempty intersection with $Q_i$, namely polygon $D^\star$, shown in fig.~\ref{fig.curv-integration-domain-convex} b). If all vertices of $P_h$ and $P_h^\star$ are regular then domain $D^\star$ consists of two triangles from $\mathcal T^Q_i$. Preimage of $D^\star$ is the quadrilateral $D$ belonging to the face $G_k$ also consisting of two triangles from $\mathcal T^G_k$. Thus
mapping $\psi_h: D^\star \to D$ is affine for each triangle from these pairs.

In the case of weak regularity one can construct more general piecewise-linear homeomorphism$\psi_h$. On fig.~\ref{fig.curv-integration-domain-concave} it is shown fragment of saddle surface. In this case dual face $Q_i$ is nonconvex and it is too restrictive to require that $Q_i$ is star-shaped with respect to $p^l_i$.

\begin{myfigure}{fig.curv-integration-domain-concave}
{
a) \hspace{4cm} b)  \\
\mbox{\includegraphics[scale=0.5, angle=0]{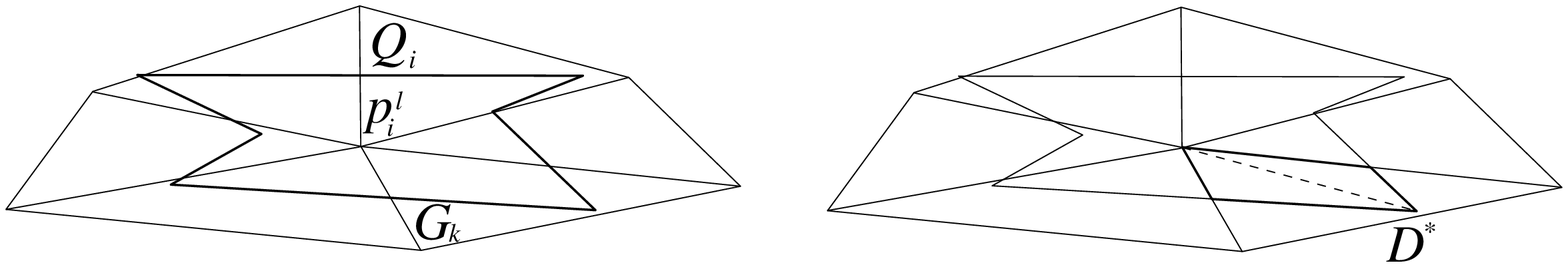} } }{Construction
of piecewise-affine homeomorphism $\psi_h: P_h^\star \to P_h$, non-convex case.}
\end{myfigure}
In this more general case dual faces can be triangulated in such a way that all angles are bounded from below and construction of homeomorphism should be based on projection and intersection of general triangulations.

Denote by $\phi_h$ and $\phi_h^\star$ piecewise affine mappings coinciding with $phi_k$ on $G_k$ and with $\phi_i$ on $Q_i$, respectively.  We will use notations $\nabla \phi_h$ and $\nabla \phi_h^\star$ for piecewise constant functions which coincide with jacobian matrices of mappings $\phi_h$ and $\phi_h^\star$ where the jacobians are well defined.

As a result the curvature deviation measure can be defined as follows
\[
\delta(p) = ||\nabla \phi_h(\psi_h(p) ) - \nabla \phi^\star_h(p) ||,
\]
where $p$ is a point lying on $P_h^\star$, and function $\delta(p)$ is piecewise constant. Here $|| \cdot||$ mean Frobenius matrix norm.

\begin{theorem}  \label{th.curv-2d-approximation-convergence}
Consider closed polyhedral surface $P_h$ which is inscribed into regular closed surface $M$ and is normal graph over $M$ realized via homeorphism $\psi$. Suppose that all faces of $P_h$ are triangles with minimal angle bounded from below, while length of edges satisfies quasi-uniformity condition
\[
 C h \leq l_j \leq h, \  0 < C < 1
\]
Let closed polyhedral surface $P_h^\star$ be normal graph over $M$ realized via homeomorphism
$\psi^\star$, and vertices of $P_h$ lie on the faces of $P_h^\star$. Suppose that all vertices of
$P_h$ and $P_h^\star$ are regular and angles of triangulations $\mathcal T^Q_i$ and $\mathcal
T^G_k$ are bounded from below, and the following closeness condition for discrete curvatures holds
\[
\sup \delta \leq O(h)
\]
i.e. for all $k$ and for every triangle $T_{km} \in \mathcal T^G_k$ the following inequality holds
\[
|| A_{ij} \nabla \psi_h - A^\star_{km}|| \leq O(h),
\]
where $T_{km} = \psi_h(T_{ij})$, $T_{ij} \in \mathcal T^Q_i$.
Then with $h \to 0$ mappings $\psi$ and $\psi^\star$ converge to identity, and
\[
|| A(p^l_i) - A^\star_{ij}|| \leq  O(h),
\]
\[
|| A({\psi^\star}^{-1}(q^l_k)) - A_{km}|| \leq O(h)
\]
\end{theorem}

One can consider average deviation of discrete measures, namely integrals
\begin{equation} \label{eq.curv-primal-dual-difference-2d-2}
\delta_2(P_h, P_h^\star) = \int \limits_{P_h} ||\nabla \phi_h(\psi(p)) - \nabla \phi^\star_h(p) ||^2 \, d \sigma
\end{equation}
and
\begin{eqnarray} \label{eq.curv-primal-dual-difference-2d-1}
\delta_1(P_h, P_h^\star) = \int \limits_{P_h} (||\nabla \phi_h(\psi(p)) - \nabla \phi^\star_h(p) ||   + \\ \nonumber
+ |\det \nabla \phi_h(\psi(p)) - \det \nabla \phi^\star_h(p) |) \, d \sigma,
\end{eqnarray}
where $d \sigma$ denote surface differential on $P_h$.

If conditions of theorem \ref{th.curv-2d-approximation-convergence} hold, but instead of maximum norm deviation one imposes closeness in a weak sense, i.e.
\[
\delta_2(P_h, P^\star_h) \to 0 \mbox{ when } h \to 0,
\]
then one can conjecture that the following weak convergence estimate holds
\begin{eqnarray} \label{eq.curv-primal-dual-convergence-2d-2}
\int \limits_{P_h} ||A(\psi^{-1}(p)) + \nabla \phi_h(p) ||^2 \, d \sigma \to 0,  \\ \nonumber
\int \limits_{P_h^\star}| |A({\psi^\star}^{-1}(p)) + \nabla \phi^\star_h(p) |^2 \, d \sigma \to 0
\end{eqnarray}

If one assumes
\[
\delta_1(P_h, P^\star_h) \to 0 \mbox{ when } h \to 0,
\]
then one can expect the following convergence estimate
\begin{eqnarray} \label{eq.curv-primal-dual-convergence-2d-1}
\int \limits_{P_h} (||A(\psi^{-1}(p)) + \nabla \phi_h(p) || +  \\ \nonumber
+ |\det A({\psi}^{-1}(p)) - \det \nabla \phi_h(p) |
)  \, d \sigma \to 0
\end{eqnarray}
and
\begin{eqnarray} \label{eq.curv-primal-dual-convergence-2d-1-dual}
\int \limits_{P_h^\star} (||A({\psi^\star}^{-1}(p)) - \nabla \phi^\star_h(p) || + \\ \nonumber
+ |\det A({\psi^\star}^{-1}(p)) - \det \nabla \phi^\star_h(p) |
)  \, d \sigma \to 0
\end{eqnarray}

Thus one can expect that the sequence of polyhedral surfaces where discretized functional (\ref{eq.curv-mean+absolute}) is bounded remains in the class of surfaces of bounded curvature.

Nice property of resulting discrete energy is that it provides exact values of curvatures for polyhedron, inscribed into sphere.
Consider polyhedra inscribed and  circumscribed around unit $d$-dimensional sphere in $\R^{d+1}$. On fig.~\ref{fig.curv-exact-circle} it is shown simplest case of a unit circle ($d=1$).

\begin{myfigure}{fig.curv-exact-circle}
{
\includegraphics[scale=0.25, angle=0]{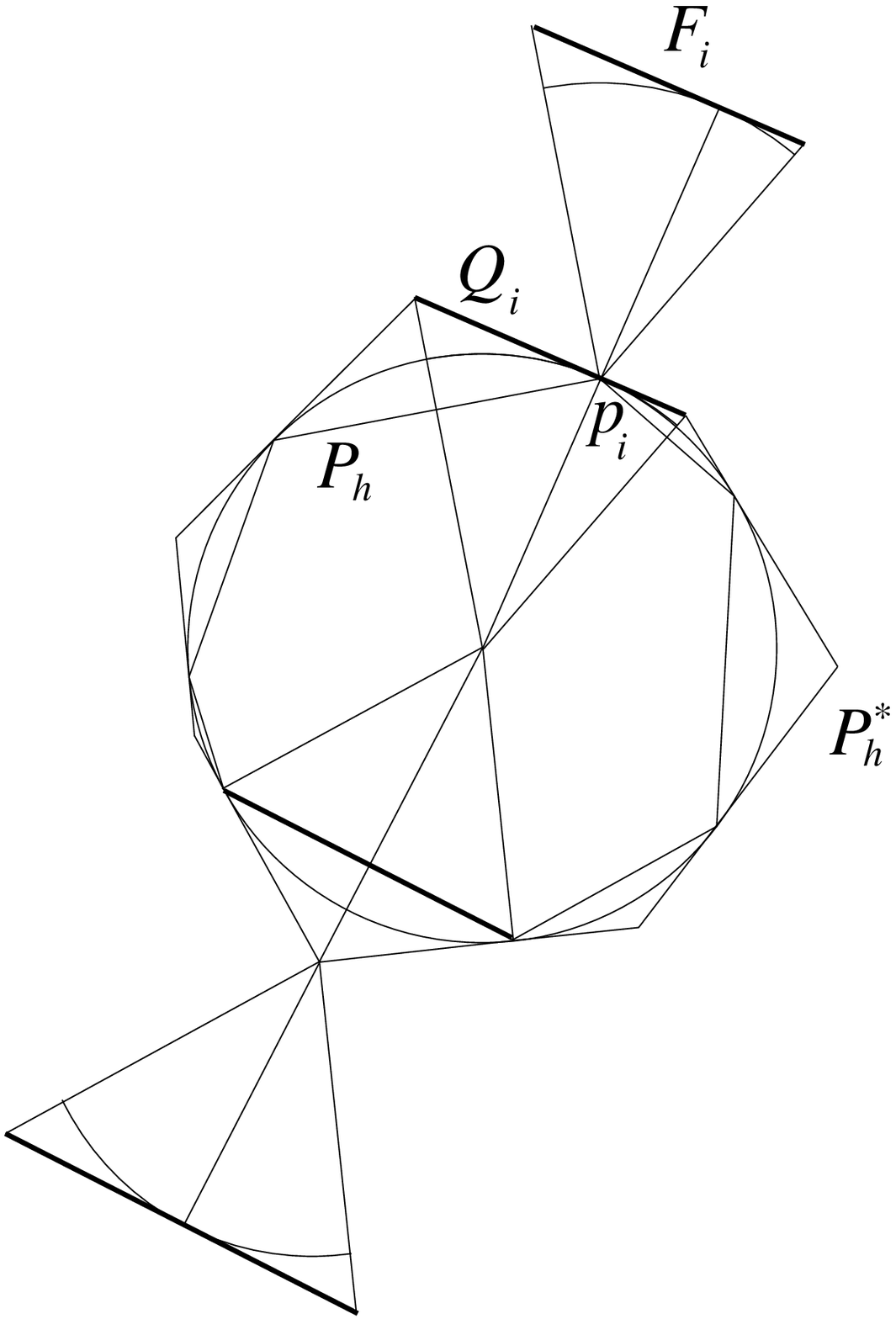}
}{For sphere dual face and normal image are congruent.}
\end{myfigure}
One can easily see that face $Q_i$, dual to vertex $p_i$ of inscribed polyhedron $P_h$ is congruent to the normal image of vertex $p_i$. Hence shape operator matrix is equal to minus identity matrix. It is also obvious that congruence property holds in $d$-dimensional case as well.

In 2004 A.I. Bobenko~\cite{Bobenko-2005-1}  introduced discrete Willmore energy $W(P_h)$ (called also conformal energy) being invariant to 3D Mebius transforms. This discrete energy is supposed to approximate exact Willmore energy
\[
W(M) = \frac12 \int\limits_M (k_1 - k_2)^2 \, d \sigma = \int\limits_M (\frac12 \tr (A^T A) - \det A) \, d \sigma,
\]
where $A$ is gradient of spherical map, i.e. shape operator matrix.

In~\cite{Bobenko-2005-1} it was proven, that $W(P_h) = 0$, when $P_h$ is convex polyhedron inscribed into sphere.

The same property holds for duality-based energy
\begin{equation} \label{eq.curv-duality-Willmore}
E(P_h^\star) = \sum\limits_i (\frac12 \tr (A^T A) - \det A)|_{Q_i} \area(Q_i)
\end{equation}
since matrix $A$ on each  dual face $Q_i$ is equal to $-I$.

Main differences with  conformal energy are that the discrete duality-based energy properties are the same for $d$-dimensional sphere in $\R^{d+1}$, and terms in discrete Willmore energy (\ref{eq.curv-duality-Willmore}) converge to the same terms in exact Willmore energy without any assumptions on polyhedra $P_h$ and $P_h^\star$ beside those guaranteeing area convergence. This is not the case for conformal energy which converges only on special polyhedra.

Duality principle allows to discern polyhedral approximations from limiting polyhedral surfaces. Let us explain this statement using simple example.
\begin{myfigure}{fig.curv-cubes}
{
a) \hspace{4cm} b)  \\
\mbox{\includegraphics[scale=0.5, angle=0]{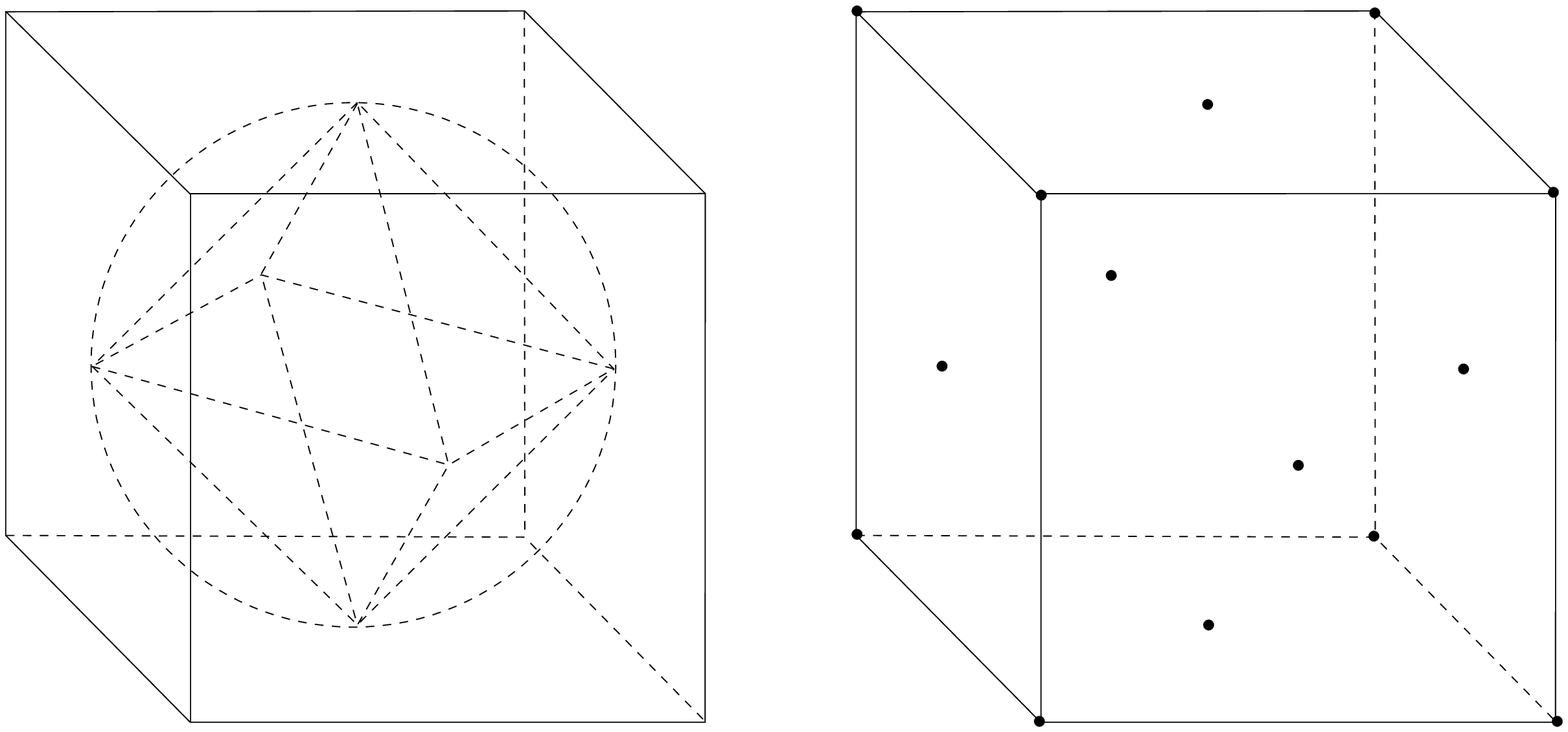}
}
}{a) Dual regular polyhedra, b) coinciding polyhedra.}
\end{myfigure}

Fig.\ref{fig.curv-cubes} a) shows cube and octahedron, dual(polar) to each other with respect to sphere. Thus one can assume that limiting surface $M$ is sphere or another regular surface. Suppose now that $P_h$ is obtained by adding vertices at the centroids of cube faces which is shown in fig.~\ref{fig.curv-cubes} b). Cube faces are triangulated including centroids. Formal construction of dual surface $P_h^\star$ leads to conclusion that it coincides with $P_h$ thus making limiting surface $M$ undistinguishable from $P_h$ and $P_h^\star$. One can compute exact value of integral (\ref{eq.curv-mean+absolute}).
\begin{equation} \label{eq.curv-edge-mean+absolute}
 E_1(P_h = P_h^\star)  = \sum\limits_e |e| 2 \tg (\frac{\varphi_e}{2}) +
\sum\limits_i \area(F_i) ,
\end{equation}
where first sum is the sum of discrete mean curvatures concentrated at edges, $\varphi_e$ denote angle between normals to faces adjacent to edge $e$, while second sum is taken over vertices.

It is clear that one can construct discrete curvature functionals using spherical image of a neighborhood of a vertex instead of normal image. As a result functional (\ref{eq.curv-mean+absolute}) can be written in the way independent on auxiliary values, such as normals at the vertices of polyhedron $M$, which define polygons $F_i$
\begin{equation} \label{eq.curv-edge-mean+absolute-spherical}
 E_1(P_h = P_h^\star)  = \sum\limits_e |e| \varphi_e +
\sum\limits_i |K_i|,
\end{equation}
where $K_i$ - absolute extrinsic (Gauss) curvature at $i$-th vertex of polyhedron.

Thus using spherical image has some advantages but computation of discrete curvatures requires complicated formulas of spherical trigonometry , while discrete curvature functional based on normal image need only computation of minors of small matrices. Note that functionals (\ref{eq.curv-edge-mean+absolute}) always majorates (\ref{eq.curv-edge-mean+absolute-spherical}). This relation is true for poledral surfaces approximating regular ones as well.

Let us remark that class of surfaces with bounded curvature measures in the sense of normal image is the subset of lipschitz continuous surfaces. It should be noted as well that method for computing discrete curvature described above can be applied in the case of surfaces with boundary. Moreover, the same duality principle can be applied for $d$-dimensional surfaces in $\R^{d+1}$ and allows to introduce extrinsic discrete curvatures in this case. In particular, the proofs of basic theorems \ref{th.curv-affine-map-curvature} and \ref{th.curv-affine-map-curvature-dual}  does not use the fact that $d=2$.

\begin{myfigure}{fig.curv-bad-parabola-tri-dual}
{
a) \hspace{3cm} b)  \\
\mbox{\includegraphics[scale=0.3, angle=0]{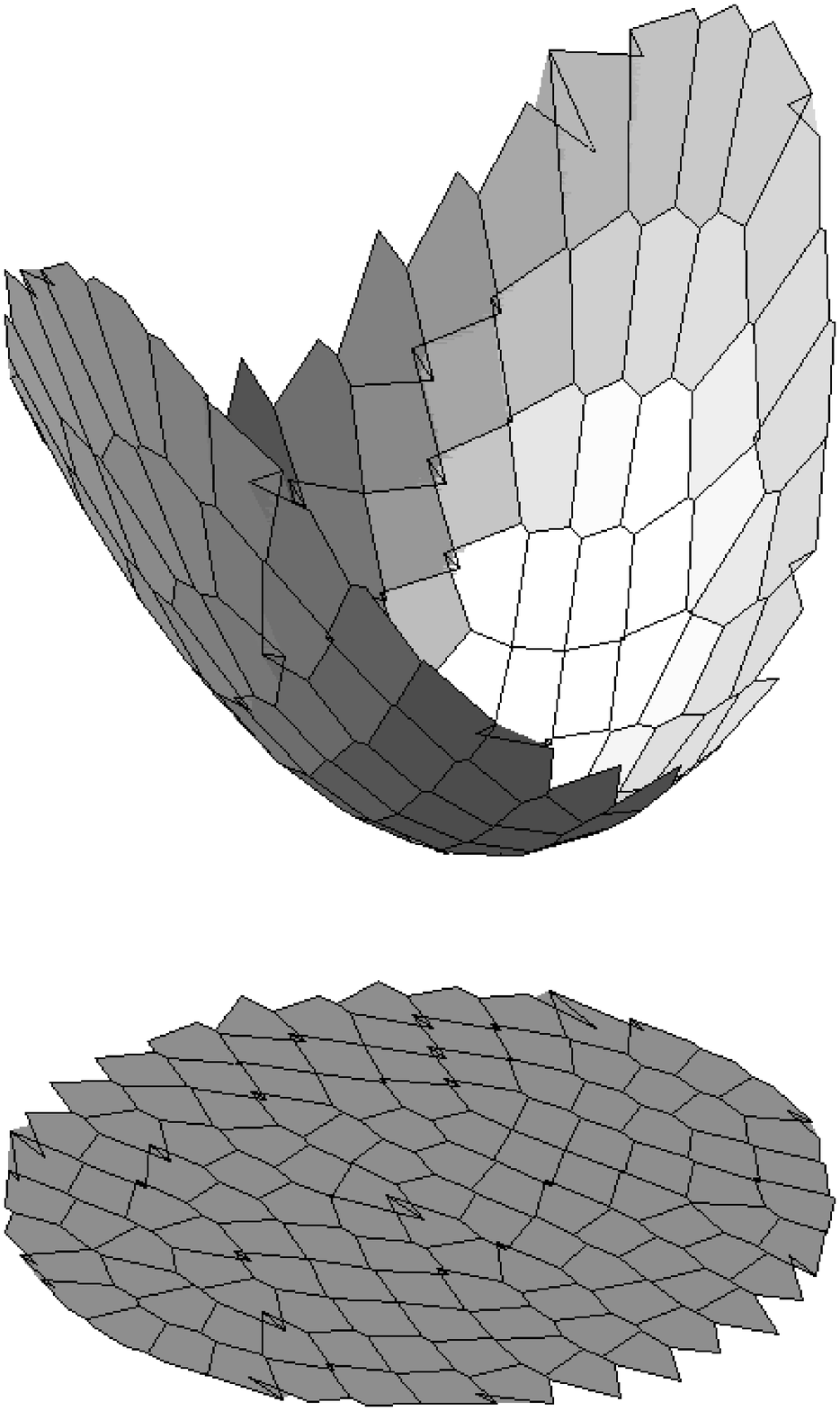}
\includegraphics[scale=0.2, angle=0]{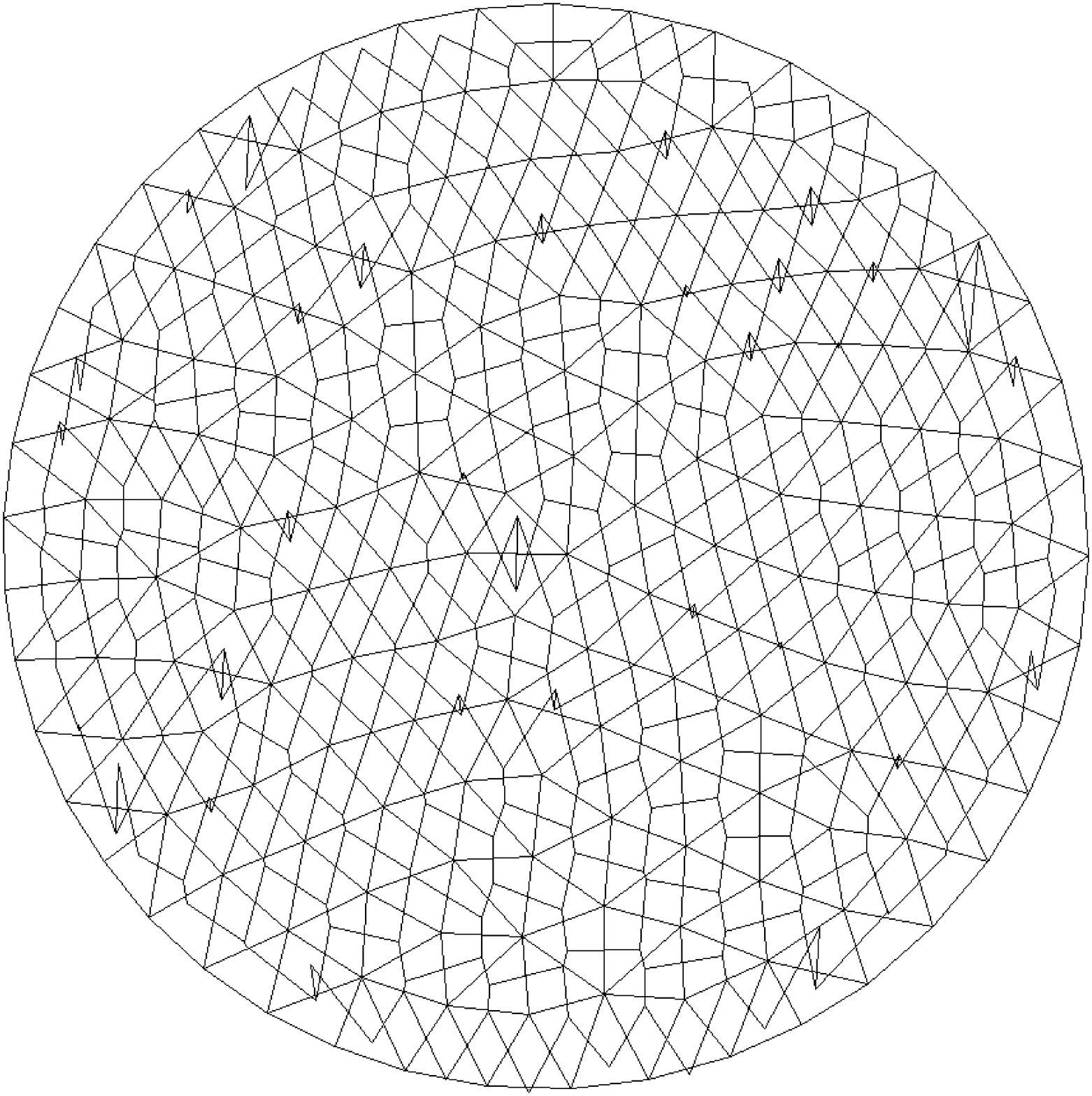}}
}{ a) Dual polyhedron with errors, b) projection of surface triangulation and dual surface.}
\end{myfigure}

{\bf Dicsrete curvature measures for non-regular polyhedra. } Polyhedron inscribed into regular
surface can contain non-regular vertices even in the case when face normals converge to exact
normals of the surface with refinement of polyhedral surface. An example of non-regular
triangulation $P_h$ inscribed into elliptic paraboloid
\[
x_3 = 2 x_1^2 + \frac15 x_2^2,
\]
is shown in fig.~\ref{fig.curv-bad-parabola-tri-dual}. One can see from the figure~\ref{fig.curv-bad-parabola-tri-dual} a) that certain ``faces'' of the dual surface $P_h^\star$ are selfintersecting, while triangulation $P_h$ consists of well shaped triangles.

Projections of surface triangulation $P_h$ and dual surface $P_h^\star$ on the plane $x_3=0$ are
shown on fig.~\ref{fig.curv-bad-parabola-tri-dual} b). Here direction of maximal curvature is
horizontal. It is clear that non-convex edges of $P_h$ lead to dual edges with wrong orientation
and hence to self-intersection. Thus eliminating ``non-convex'' edges and creation of another
edges in resulting quadrilaterals will make all vertices of  $P_h$ regular, even though may lead
to triangles with small angles.

\begin{myfigure}{fig.curv-lantern}
{
\includegraphics[scale=1.5, angle=0]{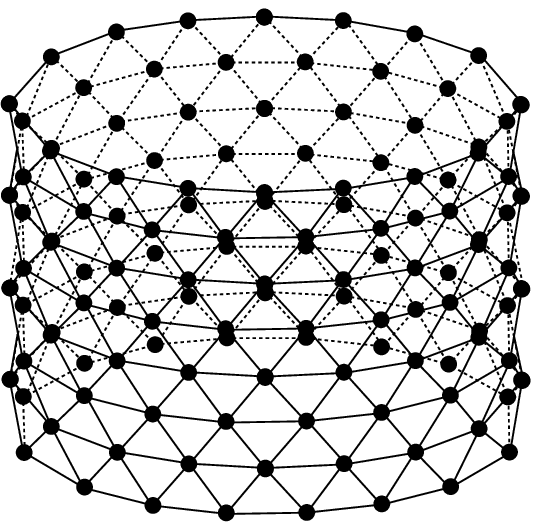}
}{Schwarz lantern: triangulation of the cylindrical surface.}
\end{myfigure}

Fig.~\ref{fig.curv-lantern} shows the so-called Schwarz lantern, i.e. polyhedron inscribed into
circular cylinder. This lantern can be constructed by subdividing the side surface of the cylinder
by planes orthogonal to its axis into $m$ equal parts. Into circle which is obtained in each cross
section one inscribes regular $n$-sided polygon. This polygon is rotated by $\pi/n$ when passing
to the next cross-section. As a result the side surface of discretized cylinder consists of
isosceles triangles. When $m,n \to \infty$ this triangulation converge to the surface of cylinder
pointwisely, but the limit of the sum of the triangle areas is $2 \pi R (H^2 + \frac{\kappa}{4}
\pi^4 R^2)^\frac12$, where $\kappa$ is the limit of the ratio $n/m^2$ when $n \to \infty$,
provided that it exists\cite{Dubnov-1961}. Here $R,H$ denote radius and height of the cylinder.

The neighborhood of nonregular vertex $p_i^l$ of a surface triangulation can be described as a
``fan'', i.e. as a cone $K^+$ with wrinkles. Normal image $\Sigma^+$ of the surface of this fan is
non-simplyconnected domain with  self-intersecting booundary, which is shown in
fig.~\ref{fig.curv-smoothed-cones}~a). If cone $K^+$ belongs to a certain half-space, one can
construct its convex envelope - cone $K^+_p$. Normal image of convex cone $K^+_p$ is convex
polygon $\Sigma^+_p$ drawn in fig.~\ref{fig.curv-smoothed-cones} a) by bold lines. We will call
$\Sigma^+_p$ by principal component of the normal image. Principal component can be constructed
for saddle point as well, which is shown in fig.~\ref{fig.curv-smoothed-cones} b).

In order to construct principal component one need to eliminate edges from $K^-$ until one obtains
canonical saddle $K^-_p$ with a normal image being quadrilateral with the sides of nonpositive
turn.

Generic cone with undefined normal image is shown in fig.~~\ref{fig.curv-smoothed-cones} c).

Extraction of principal component essentially means that discrete curvature is constructed for
some other polyhedral surface with the same set of vertices. Thus, even for approximation of
regular surface, part of the edges only prevent ``good'' approximation.

It is natural to call principal component of discrete curvature measure the discrete curvature
measure based upon principal component of the normal image.  It should be noted that principal
component of the discrete curvature measure for Schwarz lantern converges to curvature measure of
exact cylinder even though area converges to wrong value.


\begin{myfigure}{fig.curv-smoothed-cones}
{
 a) \hspace{3cm} b) \hspace{3cm} c) \\
\includegraphics[scale=0.6, angle=0]{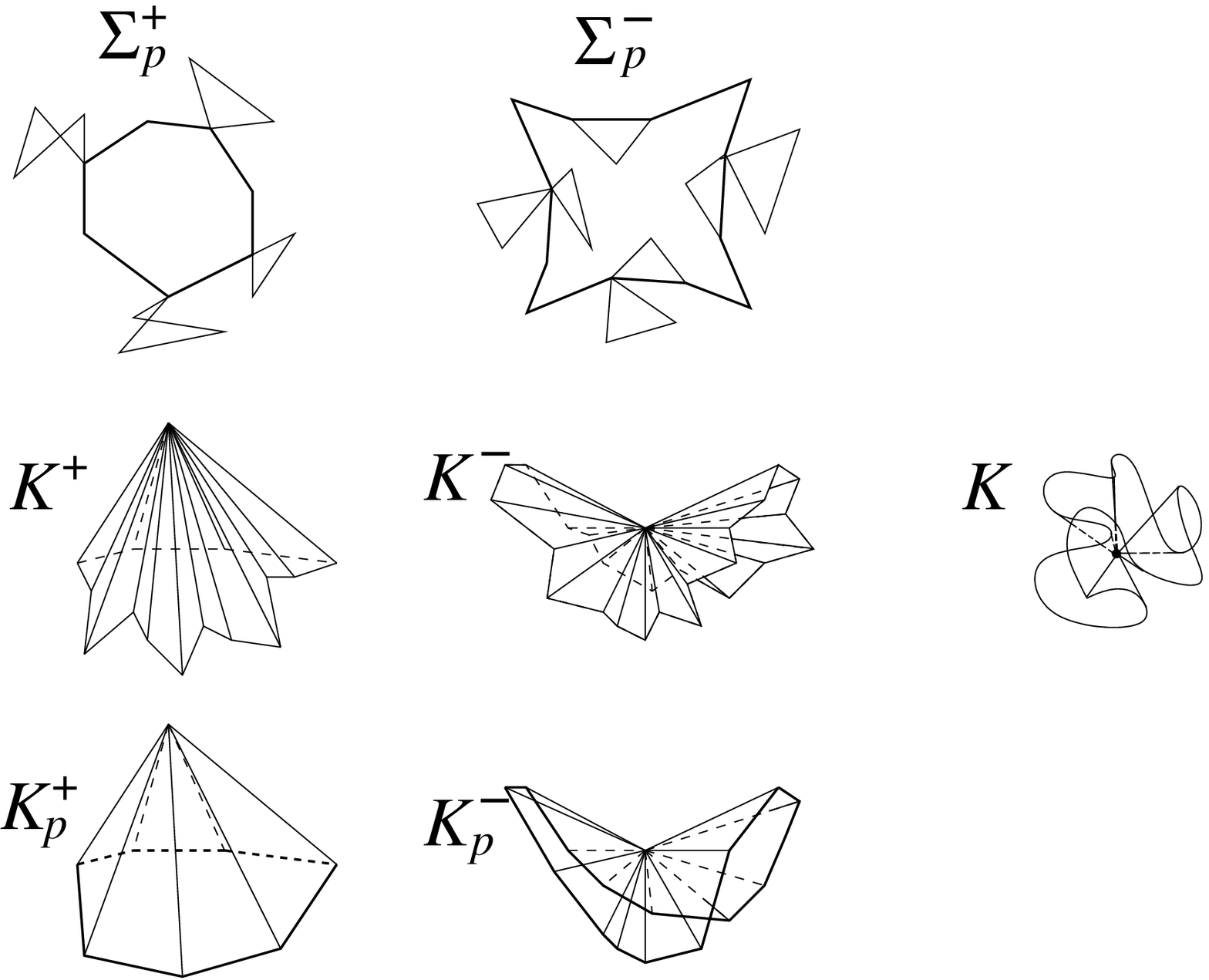}
}{ Filtering of edges and principal component of normal image.}
\end{myfigure}

When normal image $F_i$ is self-intersecting, its boundary can be decomposed into simple closed
arcs bounding simple polygons with different orientations. Summing up the areas of these polygons
with the signs corresponding to orientations of their boundary contours, one can obtain a value
which converges to $2 \pi - \sum \theta_{ij}$ with diameter of $F_i$ tending to zero. Here
$\theta_{ij}$ denote angles incident to the $i$-th vertex, and $2 \pi - \sum \theta_{ij}$ is
intrinsic (Gauss) curvature at the $i$-th vertex. One can compute the sum of absolutes values of
areas of simple polygons thus obtaining the value
\[
\area{{F_i}^p} + \delta (F_i), \ \delta(F_i) > 0,
\]
where ${F_i}^p$ denotes principal component of the normal image. In order to take into account the
original shape of polyhedral surface, one can augment principal component of the integral
curvature measure, satisfying inequality (\ref{eq.curv-generic-curvature-measure}) by a term
\[
\sum\limits_i \delta(F_i),
\]
where the sum is taken over all nonregular vertices.

It may happen that after certain edges are eliminated by filtering procedure, the cones whose
summits are adjacent vertices of polyhedral surface, become inconsistent. In this case one assumes
that certain faces of the surface $P_h$ coincide with faces of dual surface $P_h^\star$. In this
case normal images and dual faces for these singular vertices should be decomposed into
subdomains. Some sub-domains should be eliminated while discrete curvatures defined on other
subdomains allow to approximate one-sided limits of exact curvatures.

Simply put, one can formulate the optimality criteria for polyhedral approximations as follows:
optimal polyhedral surface should minimize difference between intrinsic absolute curvature and
extrinsic absolute curvature. Since extrinsic absolute curvature cannot be less than intrinsic
one, then to some extent this optimality means minimization of absolute extrinsic curvature but
unlike \cite{Alboul-2002}, this minimization is applied only at the non-regular vertices of
polyhedron. When polyhedral surface is regular in the sense of definitions \ref{def-regular-1} and
\ref{def-regular-2} its intrinsic and extrinsic curvatures coincide.

One can conclude that duality principle allows the possibility to investigate convergence of
polyhedral approximations of non-regular surfaces. It can be applied in the general case of
$d$-dimensional surface in $\R^{d+1}$.

\end{document}